\documentclass[11pt, english]{article}
\usepackage{babel}
\usepackage[latin1]{inputenc}
\usepackage[T1]{fontenc}
\usepackage{amsmath}
\usepackage{graphicx}
\usepackage{amssymb}
\usepackage{amsthm}
\usepackage{subcaption}
\usepackage{bbm}
\usepackage{mathtools}

\usepackage{xr}

\usepackage{algorithmic}
\usepackage{courier} 
\usepackage{algorithm}
\graphicspath{{figures/}}

\usepackage{chngcntr}
\usepackage{apptools}
\AtAppendix{\counterwithin{thm}{section}}

\def\^{\widehat}

\newcommand{\norm}[1]{\Vert #1 \Vert}

\def\phi{\varphi}

\numberwithin{equation}{section}

\renewcommand{\phi}{\varphi}

\def\~{\widetilde}
\def\^{\widehat}
\newcommand{\ee}{{\rm e}\hspace{1pt}}
\newcommand{\ii}{\text{i}\hspace{1pt}}
\newcommand{\dd}{\hspace{1pt}{\rm d}\hspace{0.5pt}}
\newcommand{\abs}[1]{\left| #1 \right|}

\newcommand{\veps}{\varepsilon}

\newcommand{\normm}[1]{\left\lVert#1\right\rVert}
\def\smatrix#1{\left[\begin{smallmatrix}#1\end{smallmatrix}\right]}

\newtheorem{thm}{Theorem}

\newtheorem{lem}[thm]{Lemma}
\newtheorem{cor}[thm]{Corollary}
\newtheorem{defn}[thm]{Definition}

 \newtheorem{innercustomthm}{Theorem}
\newenvironment{customthm}[1]
  {\renewcommand\theinnercustomthm{#1}\innercustomthm}
  {\endinnercustomthm}

  \newtheorem{innercustomlemma}{Lemma}
  \newenvironment{customlemma}[1]
  {\renewcommand\theinnercustomlemma{#1}\innercustomlemma}
  {\endinnercustomlemma}

\title{Sampling of Stochastic Differential Equations using the Karhunen--Lo\'{e}ve Expansion
   and Matrix Functions}

\author{Antti Koskela$^{1}$ and Samuel D. Relton$^2$ \vspace{5mm} \\ 
$^1$ Helsinki Institute for Information Technology HIIT,\\
    Department of Computer Science, University of Helsinki, Finland \\
  $^2$ 
    Leeds Institute of Health Sciences, University of Leeds, UK}

    \date{}

\begin{document}
	
\maketitle

\abstract{
   We consider linearizations of stochastic differential equations with additive noise using the Karhunen--Lo\'eve expansion.
   We obtain our linearizations by truncating the expansion and writing the solution as a series of matrix-vector products using the theory
   of matrix functions. Moreover, we restate the solution as the solution of a system of linear differential equations.
   We obtain strong and weak error bounds for the truncation procedure and show that, under suitable conditions,
   the mean square error has order of convergence $\mathcal{O}(\frac{1}{m})$ and the second moment has a
   weak order of convergence $\mathcal{O}(\frac{1}{m})$, where $m$ denotes the size of the expansion. We also discuss efficient numerical linear algebraic
   techniques to approximate the series of matrix functions and the linearized system of differential equations.
   These theoretical results are supported by experiments showing the effectiveness of our algorithms when
   compared to standard methods such as the Euler--Maruyama scheme.
}

%%%%%%%%%%%%%%%%%%%%%%%%%%%%%%
\section{Introduction}
\label{sec.introduction}
%%%%%%%%%%%%%%%%%%%%%%%%%%%%%%
Many applications of machine learning are in domains
that can be modelled using a stochastic differential equation
(SDE) with additive noise.
The diverse array of examples across many domains
includes
epidemiology~\cite{gghmp11},
weather forecasting~\cite{immm15},
finance~\cite{jack02}, and
gene expression~\cite{scra17}.
In the past decade,
advances in computing power have led to renewed
interest in this topic from the
machine learning and uncertainty quantification
communities~\cite{saso19}.

There are two common tasks when using such models.
First, the forwards propagation of uncertainty
from the parameters of the SDE
into the output allows one to calculate the
expected solution of a system,
its variance,
confidence intervals around the solution,
or indeed its entire probability distribution~\cite{Lord}.
Second, estimating the most likely parameters of a
hypothesized SDE from observations using
Bayesian inference is becoming increasingly common~\cite{bro18} .
Both tasks rely heavily upon the ability to efficiently sample
different realizations from these SDEs and typically use
Monte Carlo approaches
(e.g.~\cite{Jentzen},~\cite{Kloeden_exp},~\cite{pprs12}).

The SDEs underlying many of these applications is a
multidimensional Ornstein-Uhlenbeck process.
For example a particle of dust floating in gas
exhibits Brownian motion,
but if there is an additional convective force
generated by airflow then this Ornstein-Uhlenbeck
process can be described by the linear SDE
\begin{equation}
  \label{eq.sde_example}
  \dot{u}(t) = Lu(t) + BdW_t,
\end{equation}
where $u(t)$ is the position of the particle,
$L$ describes the convection,
and $BdW_t$ describes the stochasticity.
%We may wish to find the expected trajectory of the particle
%along with a 95\% confidence interval that bounds the region
%where we expect the particle to travel, for example.
Multidimensional Ornstein--Uhlenbeck processes of
the form~\eqref{eq.sde_example} that arise from
parabolic partial differential equations
have also been considered~\cite{Larsson}.

The primary goal of this research is to increase the speed with
which samples can be generated from
multidimensional Ornstein-Uhlenbeck processes
by exploiting a connection with matrix functions.
Although we focus on the Ornstein-Uhlenbeck process,
our approach can be extended to other stochastic processes.

A matrix function is an operator
$f \colon \mathbb{C}^{n \times n} \rightarrow \mathbb{C}^{n \times n}$
which generalizes useful properties of its scalar equivalent.
For example,
the matrix exponential is
\begin{equation*}
  \exp(A) = \sum_{k = 0}^\infty \frac{A^k}{k!}.
\end{equation*}
This generalizes many of the useful properties of the scalar exponential, for example $\exp(2A) = \exp(A)^2$.
This property forms the basis of the popular scaling and squaring method used to compute the exponential in practice~\cite{alhi11}.

The remainder of this work is organized as follows.
In section~\ref{sec.deterministic} we introduce our
approach on a deterministic semilinear ODE and derive
results that are used throughout the
rest of the analysis.
Next, in section~\ref{sec.karhunen-loeve},
we apply our approach to a
SDE with additive noise.
We define our numerical method in section~\ref{sec.numerical-scheme}
and derive
error bounds on the expected truncation error introduced.
In section~\ref{sec.matexp} we reformulate our numerical method as a
single matrix--vector product involving the matrix exponential
and implementation strategies.
Numerical experiments are given in
section~\ref{sec.numerical-experiments}
with concluding remarks in section~\ref{sec:conclusion}.

%%%%%%%%%%%%%%%%%%%%%%%%%%%%%
\section{Deterministic differential equations}
\label{sec.deterministic}
%%%%%%%%%%%%%%%%%%%%%%%%%%%%%
We begin by deriving our approach for a
deterministic semilinear differential equation before introducing
the additional complexity of the stochastic terms.
In brief,
we want to relate the solution of the semilinear differential equation
\begin{equation} \label{eq:semilinearODE}
  \dot{u}(t) = L u(t) + g(t), \, u(0) = u_0 \in \mathbb{R}^n,
  \, L \in \mathbb{R}^{n \times n},
\end{equation}
to matrix functions when the nonlinear function $g(t)$ is
approximated by a finite dimensional Fourier series.
The approach is reminiscent
of~\cite[Thm.\;2.1]{alhi11}. %~and~\cite[Lemma\;1]{koskela2013}.
To this end, we consider the Fourier series of the nonlinear function
$g \, : \, \mathbb{R} \rightarrow \mathbb{R}^n$,
\begin{equation*}
g(t) = \sum\nolimits_{k=0}^\infty a_k \cos( c_k t) + b_k \sin ( c_k t),
\end{equation*}
where $a_k,b_k \in \mathbb{R}^n$ and $c_k \in \mathbb{R}$, $k\geq 1$.

%%%%%%%%%%%%%%%%%%%%%%%%%%%%%%
\subsection{Truncated system and matrix functions}
\label{sec.varphi-trig}
%%%%%%%%%%%%%%%%%%%%%%%%%%%%%%
Now we can replace $g(t)$ in the
semilinear equation~\eqref{eq:semilinearODE}
by a truncated Fourier series expansion
\begin{equation} \label{eq:truncated}
  g(t) \approx g_N(t) = \sum\nolimits_{k=0}^N a_k \cos( c_k t) + b_k \sin( c_k t)
\end{equation}
to obtain the approximative ODE
\begin{equation} \label{eq:semilinearODE_N}
\dot{u}_N(t) = L u_N(t) + g_N(t), \quad u_N(0) = u_0 \in \mathbb{R}^n.
\end{equation}
The error introduced by this truncation will be bounded in the
proceeding subsection.
Using the variation-of-constants formula,
the solution of \eqref{eq:semilinearODE_N} is
\begin{equation} \label{eq:phi_series}
	\begin{aligned}
		u_N(t) &= \ee^{t L} u_0 + \sum\nolimits_{k=1}^N  \left( \int\nolimits_0^t \ee^{(t-s)L} \cos( c_k s) \, \dd s \right)  a_k \\
		 & \phantom{= \ee^{t L} u_0  \sum\nolimits_{k=1}^N\,} + \left( \int\nolimits_0^t \ee^{(t-s)L} \sin(c_k s) \, \dd s \right) b_k \\
		   &= \ee^{t L} u_0 + \sum\nolimits_{k=1}^N   \varphi_{k,t}^{\cos}(L) \, a_k +   \varphi_{k,t}^{\sin}(L) \, b_k,
	\end{aligned}
\end{equation}
	where $\varphi_{k,t}^{\cos}$ and $\varphi_{k,t}^{\sin}$ denote the  functions
\begin{equation*} % \label{eq.phi_integral_defn}
	\begin{aligned}
	\varphi_{k,t}^{\cos}(z) &= \int\nolimits_0^t \ee^{(t-s)z} \cos(c_k s) \, \dd s, \\
    \varphi_{k,t}^{\sin}(z) &= \int\nolimits_0^t \ee^{(t-s) z} \sin(c_k s) \, \dd s.
	\end{aligned}
\end{equation*}
These functions are clearly analytic on the whole complex plane
and for $z \in \mathbb{C}$ satisfy
\begin{equation} \label{eq:phi}
	\begin{aligned}
		\varphi_{k,t}^{\cos}(z) + \ii \varphi_{k,t}^{\sin}(z) & = \int\nolimits_0^t \ee^{(t-s)z} \ee^{\ii c_k s } \, \dd s \\
		 &= \frac{\ee^{t z }- \ee^{ \ii c_k t}}{z-\ii c_k}.
	\end{aligned}
\end{equation}
We call $\phi_{k,t}^{\cos}$ and $\phi_{k,t}^{\sin}$ the
\emph{trigonometric $\varphi$ functions}.
% due to the similarity of the above
% integral formulation~\eqref{eq.phi_integral_defn} and the
% $\varphi_{\ell}$ functions that occur in exponential integrators
% related to the matrix exponential~\cite{niwr12}
% \begin{equation*}
%   \varphi_{\ell}(z) =
%   \frac{1}{(\ell - 1)!} \int\limits_0^1 e^{(1-s)z}s^{\ell - 1}ds.
% \end{equation*}
If we consider the real and imaginary parts of \eqref{eq:phi} separately,
we can make the following definition of the
(scalar versions of the) trigonometric $\varphi$ functions
without an integral.
\begin{defn}
  \label{def.phitrig}
  Let $(\lambda_k)_{k=1}^\infty$ be a sequence of real numbers.
  Then the trigonometric $\varphi$ functions
  corresponding to the sequence $(\lambda_k)_{k=1}^\infty$ are defined
  for $z\in \mathbb{C}$ and $t\in \mathbb{R}$ by
\begin{equation}  \label{eq:phi_cos}
\varphi_{k,t}^{\cos}(z)  = \frac{ z \ee^{z t} - z \cos(\lambda_k t) + \lambda_k \sin(\lambda_k t)}{z^2+\lambda_k^2},
\end{equation}
and
\begin{equation} \label{eq:phi_sine}
\varphi_{k,t}^{\sin}(z)  = \frac{\lambda_k \ee^{z t}- z \sin(\lambda_k t)- \lambda_k \cos(\lambda_k t)}{z^2+\lambda_k^2}.
\end{equation}
\end{defn}
%From \eqref{eq:phi} it clearly follows that these functions are analytic on the whole complex plane.

%%%%%%%%%%%%%%%%%%%%%%%%%%%%%%
\subsection{Bound for the truncation error}
%  using $\varphi$ functions}
\label{sec.bound-trunc-error}
%%%%%%%%%%%%%%%%%%%%%%%%%%%%%%

At this point we have shown how one can form an approximation to our
original differential equation by truncating a Fourier series
expansion of the nonlinear term~\eqref{eq:semilinearODE_N}
and using matrix functions~\eqref{eq:phi_series}.
In this section we analyze the convergence of the approximate solution
$u_N(t)$ as $N$ grows.

Our analysis requires the use of the \emph{numerical range}
$\mathcal{F}(A)$ of a matrix $A \in \mathbb{C}^{n \times n}$,
which is defined as
$$
\mathcal{F}(A) = \{ x^* A x \, : \, x \in \mathbb{C}^n, \norm{x}_2 = 1 \}.
$$
We also require the related notion of the
\emph{logarithmic norm} of a matrix $A$, defined as
 $$
 \mu(A) = \max \{ \mathrm{Re} z \, : \, z \in \mathcal{F}(A) \}.
 $$
Using these, we state the following bound for the norm of the
matrix functions $\varphi_{k,t}^{\cos}(A)$ and $\varphi_{k,t}^{\sin}(A)$.

\begin{lem} \label{lem:norm_bound}
 Let $A\in \mathbb{R}^{n \times n}$, $(\lambda_k)_{k=1}^\infty \subset \mathbb{R}$ and let the corresponding
 $\varphi_{k,t}^{\cos}(z)$ and $\varphi_{k,t}^{\sin}(z)$ be defined as in \eqref{eq:phi_cos} and \eqref{eq:phi_sine}.
 Suppose $d\big(\ii \lambda_k,\mathcal{F}(A) \big) > 0$, where
 $d(c,X)$ denotes the Euclidean distance of $c\in \mathbb{C}$ from the set $X$. Then,
 $$
 \norm{\varphi_{k,t}^{\cos}( A)} \leq \frac{1 + \ee^{t \mu(A)}}{ d ( \ii \lambda_k,  \mathcal{F}(A))},
 $$
 and the same bound holds for $\norm{\varphi_{k,t}^{\sin}( A )}$.
 %where
 %If $A \in \mathbb{R}^{n \times n}$, then the above bound holds for $\norm{\varphi_k^{\cos}( t A)}$ and $\norm{\varphi_k^{\sin}( t A)}$.
 \begin{proof}
  From the representation \eqref{eq:phi} we see that
  \begin{equation} \label{eq:var_cos_A}
  \varphi_{k,t}^{\cos}( A) = \mathrm{Re} \left( (\ee^{tA} - \ee^{\ii \lambda_k t} I)(A - \ii \lambda_k I)^{-1} \right)
  \end{equation}
  and similarly $\varphi_{k,t}^{\sin}( A)$ is given by the imaginary part.
  To bound \eqref{eq:var_cos_A},
  we use the well known bounds $\norm{\ee^{A}}_2 \leq \ee^{\mu(A)}$
  (see e.g.~\cite[Sec.\;14]{Trefethen_Embree}),
  and $\norm{(zI - A)^{-1}} \leq d(z,\mathcal{F}(A))^{-1}$.
 \end{proof}
\end{lem}

To illustrate the utility of this bound,
we give a short example of how one can bound
$\norm{\varphi_{k,t}^{\cos}(A)}$%.
and a second example which applies this lemma to bound
the truncation error $\norm{u(t) - u_N(t)}$.

\textbf{Example.} Let $A \in \mathbb{R}^{n \times n}$ be a negative semi-definite matrix (i.e. $\mathcal{F}(A) \subset \mathbb{R}_{\leq 0}$, $\mu(A) \leq 0$)
and let $(\lambda_k)_{k=1}^\infty \subset \mathbb{R} \setminus \{ 0 \}.$
Then using the lemma above, we have
\begin{equation*} % \label{eq:bound_nsd}
 \norm{\varphi_{k,t}^{\cos}(A)} \leq \frac{1}{\abs{\lambda_k}}.
\end{equation*}

\textbf{Example 2.} Let $A \in \mathbb{R}^{n \times n}$ be a
negative semi-definite matrix and let
$g(t) = f(t)p$, $p \in \mathbb{R}^n$,
where $f(t)$ is a 2$\ell$-periodic "sawtooth wave", i.e.,
$$
f(t) = \frac{t}{2\ell} , \quad 0 \leq t \leq 2\ell.
$$
In this case $f(t)$ has the Fourier series defined by $b_k = -\frac{1}{k \pi}$,
$c_k = \frac{k \pi}{\ell}$.
Furthermore,
let $u(t)$ denote the solution of \eqref{eq:semilinearODE}
and $u_N(t)$ that of \eqref{eq:semilinearODE_N}.
Then using the lemma above we have
\begin{equation*} \label{eq:bound_nsd}
\begin{aligned}
& \norm{u(t) - u_N(t)}^2 \leq \sum\limits_{k = N}^\infty \frac{1}{k \pi} \norm{\varphi_k^{\sin}( A )} \norm{p} \\
 &  \quad \leq \norm{p} \ell \sum\limits_{k = N}^\infty \frac{1}{(k \pi)^2} \leq \norm{p} \ell \int\limits_{N-1}^\infty \frac{1}{(\pi x)^2} \, \dd x \\
& \quad = \frac{\norm{p} \ell} {\pi^2(N-1)}.
 \end{aligned}
\end{equation*}
% Here the smoothness result, few lines.

%%%%%%%%%%%%%%%%%%%%%%%%%%%%%%
\section{Approximation of the It$\hat{\textrm{o}}$ integral using the Karhunen--Lo\'eve expansion}
\label{sec.karhunen-loeve}
%%%%%%%%%%%%%%%%%%%%%%%%%%%%%%

In the previous section we were focused on
deterministic differential equations in order to explain our approach.
We are now ready to consider a linear differential equation with
an additive stochastic term by applying the same methodology.

Let us consider the stochastic differential equation
\begin{equation} \label{eq:semilinearSDE}
d X_t = L X_t + B dW_t,
\end{equation}
where $X_t \in \mathbb{R}^n$, $L \in \mathbb{R}^{n \times n}$ and
$W_t$ is the standard Wiener process.
The exact solution is given by (see~\cite[Sec.\;4.8]{kloeden})
\begin{equation*} % \label{eq:Ito1}
X_t = \ee^{t L} X_0 + \int_0^t \ee^{(t-s) L} B d W_s,
\end{equation*}
where the stochastic integral is defined here as the It$\hat{\textrm{o}}$ integral.
That is,
\begin{equation*} %\label{eq:Ito2}
\int _{0}^{t} f(s) \,d W_s =\lim _{n\rightarrow \infty }\sum _{[t_{i-1},t_{i}]\in \pi_{\ell}}  f(t_{i-1})(W_{t_{i}}-W_{t_{i-1}}),
\end{equation*}
where $\pi_\ell$ is a sequence of partitions of $[0,t ]$ with mesh
size going to zero as $\ell \rightarrow \infty$,
i.e.
$$
\pi_\ell = \{ [t_0, t_1], \ldots, [t_{\ell-1},t_\ell] \}
$$
such that $0=t_0 < t_1 < \ldots < t_\ell = t$.
We note that for equation \eqref{eq:semilinearSDE}
the It$\hat{\textrm{o}}$ and the Stratonovich definitions are actually equivalent.

To relate the stochastic process to the nonlinear function $g(t)$
used in the previous section,
we proceed by replacing $W_t$ with its Karhunen--Lo\'eve expansion
\begin{equation} \label{eq:series2}
W_t = \sqrt{2} \sum\nolimits_{k=1}^\infty Z_k \frac{\sin(\lambda_k t) }{ \lambda_k},
\end{equation}
where $ \lambda_k = \left(k - 1/2\right) \pi$
and $Z_k$ are independent normally distributed
vector valued random variables with zero mean and unit variance.
This also means that each component $(Z_k)_i$, $1\leq i\leq n$,
of the vectors $Z_k$ are independent
$\mathcal{N}(0,1)$ distributed random variables.

% We can now apply the approach explained in section~\ref{sec.deterministic},
% utilizing the fact that for any
% differentiable function $f$ and integrable function $g$
% the Riemann--Stieltjes integral gives
% \begin{equation} \label{eq:Stieltjes}
% \begin{aligned}
% &\int\limits_0^t g(t) df(t)\\
%  & =\lim_{n\rightarrow \infty }\sum _{[t_{i-1},t_{i}]\in \pi_{n}}  g({t_{i}})(f({t_{i}})-f({t_{i-1}})) \\
% & = \lim _{n\rightarrow \infty }\sum _{[t_{i-1},t_{i}]\in \pi_{n}}  (t_i - t_{i-1}) g({t_{i}})\frac{f({t_{i}})-f({t_{i-1}})}{t_i - t_{i-1}} \\
% & =  \int_0^t g(s) \frac{\dd }{\dd s} f(s) \, \dd s.
% \end{aligned}
% \end{equation}
To obtain a practical method we approximate the Brownian motion $W_t$
by the \emph{truncated Karhunen--Lo\'eve expansion}
\begin{equation}
  \label{eq:weiner_khl}
W_t \approx W_t^m := \sqrt{2} \sum\nolimits_{k=1}^m Z_k \frac{\sin(\lambda_k t) }{ \lambda_k}.
\end{equation}
The truncated expansion $W_t^m$ is always differentiable with
respect to $t$ and therefore%,
%in a similar way to \eqref{eq:Stieltjes},
%the truncation gives
\begin{equation} \label{eq:truncated_integral}
\begin{aligned}
& \int_0^t \ee^{(t-s) L} B \, d W_s^m \\
%& = \sqrt{2} \lim_{n\rightarrow \infty }\sum _{[t_{i-1},t_{i}]\in \pi_{n}}
%\ee^{(t-t_{i-1}) L} \left( W_{t_i}^m - W_{t_{i-1}}^m \right) \\
& = \sqrt{2} \sum\nolimits_{k=1}^m \int_0^t \left( \ee^{(t-s) L}
\left( \frac{\dd}{\dd s} \frac{\sin(\lambda_k s) }{ \lambda_k} \right) \, \dd s \right) B Z_k \\
%& = \sqrt{2} \sum_{k=1}^m \int_0^t \left( \ee^{(t-s) L} \cos( \lambda_k s) \, \dd s \right) B Z_k \\
& = \sqrt{2} \sum\nolimits_{k=1}^m  \varphi_{k,t}^{\cos}(L) B Z_k,
\end{aligned}
\end{equation}
where the functions $\varphi_{k,t}^{\cos}(z)$ correspond to the
sequence $(\lambda_k)_{k=1}^\infty \subset \mathbb{R}$
taken from the Karhunen--Lo\'{e}ve expansion.
As $m$ goes to infinity the integral converges in the $L^2$ sense,
the proof of which follows from the Wong--Zakai
theorem~\cite{stva72},~\cite{woza65a},~\cite{woza65b}.
For the convergence properties of general bases we refer to the the appendix of~\cite{Lyons}.

We are now ready to present the solution to our
stochastic differential equation in terms of
matrix functions. For the proof of the following theorem we refer to~\cite[Sec.\;2]{Twardowska}.
\begin{thm}
  \label{thm.sde_soln}
 %Suppose $L$ is negative semidefinite. Then,
 The solution of \eqref{eq:semilinearSDE} has the representation
 \begin{equation} \label{eq:series_representation}
\begin{aligned}
 X_t = \ee^{t L} X_0 + \sqrt{2} \sum\nolimits_{k=1}^\infty  \varphi_{k,t}^{\cos}(L) B Z_k,
 \end{aligned}
\end{equation}
where $Z_k$ are independent $\mathcal{N}(0,1)$ distributed vector
valued random variables, the functions $\varphi_{k,t}^{\cos}$ are
defined by the coefficients $\lambda_k$ of the Karhunen--Lo\'eve
expansion of $W_t$ and equation \eqref{eq:phi_cos}, and where the
convergence is in $L^2$ and uniform in $t$.
\end{thm}
From the representation \eqref{eq:series_representation} we may deduce the following properties.
The proof is left to the appendix.
% \begin{equation} \label{eq:exp_X_t}
% \end{equation}
% For the second moment $\mathbb{E}( \norm{X_t}^2 )$ we have the following.
\begin{thm} \label{thm:variance_thm}
We have
\begin{enumerate}%[label=(\roman*)]
\item
$
\mathbb{E}(X_t) = \ee^{t L} X_0.
$
\item
$
\mathbb{E} \left( \norm{X_t}^2 \right) = \norm{\ee^{t L} X_0}^2 + 2 \sum\limits_{k=1}^\infty \norm{ \varphi_{k,t}^{\cos}(L) B }_F^2,
$
where $\norm{ \cdot }_F$ is the Frobenius norm.
\end{enumerate}
\end{thm}

% Theorem~\ref{thm:variance_thm} implies that the trace of the covariance matrix can indeed be computed using the trigonometric $\varphi$ functions.
% \begin{cor} \label{cor:covariance}
% We have the following representation for the trace of the covariance matrix of $X_t$:
% \begin{equation}
% 	\begin{aligned}
% 		\mathrm{tr} \,\, \mathrm{Cov}(X_t,X_t) & = \mathbb{E}(\norm{X_t}^2) - \norm{\mathbb{E}(X_t)}^2 \\
% 		 &= 2 \sum\limits_{k=1}^\infty \norm{ \varphi_{k,t}^{\cos}(L) B }_F^2.
% 	\end{aligned}
% \end{equation}
% \end{cor}
For the special case when $L$ is a normal operator we have the following.
\begin{cor} \label{cor:normal}
Suppose $B = c \, I$ for some $c \in \mathbb{R}$ and that $L \in \mathbb{R}^{n \times n}$ is normal, i.e. unitarily diagonalizable.
Then,
$$
\mathbb{E} \left( \norm{X_t}^2 \right) = \norm{\ee^{t L} X_0}^2
+ 2 c^2 \sum\limits_{k=1}^\infty \sum\limits_{\lambda \in \Lambda(L)} \abs{\varphi_{k,t}^{\cos}(\lambda)}^2,
$$
where $\Lambda(L)$ denotes the spectrum of $L$.
\begin{proof}
The claim follows from Theorem~\ref{thm:variance_thm}, the fact that for all
$A \in \mathbb{R}^{n \times n}$, $\norm{A}_F^2 = \sum_i \sigma_i^2$, where $\sigma_i$'s are the singular values of
$A$, and from the fact that for normal matrices the singular values equal the absolute values of the eigenvalues.
\end{proof}
\end{cor}

%%%%%%%%%%%%%%%%%%%%%%%%%%%%%%%%%%%%%%%%%%%%%%%%%%%%%%%%%%%%%%%%%%%%%%%%%%%
\section{Numerical method and its analysis}
\label{sec.numerical-scheme}
%%%%%%%%%%%%%%%%%%%%%%%%%%%%%%%%%%%%%%%%%%%%%%%%%%%%%%%%%%%%%%%%%%%%%%%%%%%
In Theorem~\ref{thm.sde_soln} we expressed the solution
to our stochastic differential equation as
an infinite sum of matrix functions.
To evaluate this numerically we will approximate this by
a finite sum of matrix functions
 \begin{equation} \label{eq:finite_sum}
 X^m_t = \ee^{t L} X_0 + \sqrt{2} \sum\limits_{k=1}^m  \varphi_{k,t}^{\cos}(L) B Z_k.
 \end{equation}

 To simplify the discussion,
 let us assume that $L$ is negative semidefinite for the moment.
 %We will address the more general case of
 %sectorial matrices in subsection~\ref{subsec:sectorial}.
 In this case, the strong mean squared error has the following bound.
\begin{thm} \label{thm:strong_bound}
 Let the linear operator $L \in \mathbb{R}^{n \times n}$
 in \eqref{eq:semilinearSDE} be negative semidefinite.
 The error introduced by the approximation~\eqref{eq:finite_sum}
 with $m$ terms satisfies the following theorem, of which proof is in the appendix.
 \begin{equation*}% \label{eq:strong_bound}
	 \begin{aligned}
	    \mathbb{E}\left(\norm{X_t - X^m_t}^2 \right) &= 2 \sum\nolimits_{k={m+1}}^\infty \normm{\varphi_{k,t}^{\cos}(L) B }_F^2 \\
		& \leq \frac{2 \norm{B}^2}{\pi^2} \frac{n}{m-1}.
	 \end{aligned}
 \end{equation*}
\end{thm}

For the second moment of the norm of $X_t^m$,
analogously to Theorem~\ref{thm:variance_thm},
we have the following result.
\begin{lem} \label{lem:variance_thm_finite}
  $X_t^m$ satisfies the following
  \begin{equation*}
    \mathbb{E} \left( \norm{X_t^m}^2 \right) =
    \norm{\ee^{t L} X_0}^2 +
    2 \sum\nolimits_{k=1}^m \norm{ \varphi_{k,t}^{\cos}(L) B }_F^2.
  \end{equation*}
\end{lem}
It is easy to see that the weak error
$
\norm{\mathbb{E}(X_t - X^m_t)}% = \normm{\mathbb{E}\left( \sqrt{2} \sum_{k=m+1}^\infty  \varphi_{k,t}^{\cos}(L) B Z_k \right)}
$
is always zero.
For the weak error of the second moment of the solution we get the following bound,
which is a direct consequence of Theorem~\ref{thm:variance_thm} and Lemma~\ref{lem:variance_thm_finite}.
\begin{thm} \label{thm:weak_bound}
The following weak error bound holds:
\begin{equation*}
	\begin{aligned}
	    \mathbb{E}\norm{X_t}^2 - \mathbb{E}\norm{X_t^m}^2 & = 2 \sum\nolimits_{k=m+1}^\infty \norm{ \varphi_{k,t}^{\cos}(L) B }_F^2 \\
	 &    \leq \frac{2 \norm{B}^2}{\pi^2} \frac{n}{m-1}.
	\end{aligned}
\end{equation*}
\end{thm}
Theorem~\ref{thm:strong_bound} shows that the method has the same strong order of convergence
$\frac{1}{2}$ as the Euler--Maruyama method in the sense that it converges pathwise as $O(1/\sqrt{m})$ with respect to the number of time steps $m$.
Similarly, Theorem~\ref{thm:weak_bound} indicates that the second moment of the norm of the numerical solution has a weak
order of convergence $1$ (see~\cite{kloeden} for the definitions).

% Similarly to Corollary~\ref{cor:covariance}, we see that
% \begin{equation}
% 	\begin{aligned}
% 		\mathrm{tr} \,\, \mathrm{Cov}(X_t^m,X_t^m) & = \mathbb{E}(\norm{X_t^m}^2) - \norm{\mathbb{E}(X_t^m)}^2 \\
% 		& = 2 \sum\limits_{k=1}^m \norm{ \varphi_{k,t}^{\cos}(L) B }_F^2.
% 	\end{aligned}
% \end{equation}
% This implies the following weak error bound

%%%%%%%%%%%%%%%%%%%%%%%%%%%%%%
\subsection{Sectorial matrices}
\label{subsec:sectorial}
%%%%%%%%%%%%%%%%%%%%%%%%%%%%%%

The above bounds can be easily generalized to coefficient matrices $L$
which are \it sectorial. \rm This means that the numerical range of $L$
lies within a cone of a given angle opening to the left.
These matrices often occur following spatial discretizations
of parabolic PDEs that lead to nonsymmetric or
nonnormal coefficient matrices $L$,
e.g.~in advection diffusion equations.

Let $\alpha \in [0,\frac{\pi}{2})$ and define
$$
\mathcal{S}_\alpha := \{ 0 \} \cup \{ \, z \in \mathbb{C} \, : \, \abs{\mathrm{Arg}(-z) } \leq \alpha \},
$$
i.e., $\mathcal{S}_\alpha$ is a cone of angle $2 \alpha$ with its vertex in origin. Using this notation we give the
following definition.
\begin{defn}
Let $\alpha \in [0,\frac{\pi}{2})$ and $\gamma \in \mathbb{R}$. The matrix $L$ is called \rm sectorial \it
with half-angle $\alpha$ and vertex $\gamma$ if the numerical range $\mathcal{F}(L - \gamma I)$ is contained in
$S_\alpha$.
\end{defn}

If $L$ is sectorial with vertex $\gamma$ and half-angle $\alpha$,
and if $\lambda \in \mathbb{R}$ such that $\ii \lambda \not\in \mathcal{S}_\alpha$, % and $\mathcal{F}(L) \subset \mathcal{S}_\alpha$,
then by simple geometry it can be shown that
$$
\frac{1}{d(\ii \lambda, \mathcal{F}(L))} \leq \frac{1}{\abs{\lambda} \cos \alpha - \gamma \sin \alpha}.
$$

%%%%%%%%%%%%%%%%%%%%%%%%%%%%%%
\subsection{The Brownian bridge and other stochastic processes}
\label{sec.brownian-bridge}
%%%%%%%%%%%%%%%%%%%%%%%%%%%%%%

Another example of where our theory can be applied is
the Brownian bridge $B_t = W_t - t W_1$,
which can be represented as the series
$$
B_t = \sqrt{2} \sum\limits_{k=1}^\infty Z_k \frac{\sin(\lambda_k t) }{\lambda_k},
$$
where the $Z_k$ are vectors with elements drawn from
a normal $\mathcal{N}(0,1)$ distribution
and $\lambda_k = \pi k$.
Note that this is identical to the Karhunen--Lo\'{e}ve expansion
for standard Brownian motion~\eqref{eq:series2},
except that in the previous case $\lambda_k = (k-1/2)\pi$.
One can trivially adapt the results from the previous sections to
this,
and similar,
stochastic processes in order to work with a
range of different models.

%%%%%%%%%%%%%%%%%%%%%%%%%%%%%%%%%%%%%%%%%%%%%%%%%%%%%%%%%%%%%%%%%%%%%%%%%%%
\subsection{Gaussian Processes}
\label{sec.lyapunov}
%%%%%%%%%%%%%%%%%%%%%%%%%%%%%%%%%%%%%%%%%%%%%%%%%%%%%%%%%%%%%%%%%%%%%%%%%%%

We remark that the solution of the linear SDE is a Gaussian process so that it is
uniquely determined by its mean and covariance which satisfy a vector valued linear differential equation and a
matrix valued Lyapunov differential equation, respectively~\cite[Sec.\;6]{saso19}.
Thus, our proposed sampling method can be seen as a way to avoid the expensive solving of the matrix valued
differential equations. Moreover, changing the initial value of the system affects only the first term $\ee^{t L} X_0$
which allows efficient sampling also in the case the initial value $X_0$ is a random variable.

%%%%%%%%%%%%%%%%%%%%%%%%%%%%%%%%%%%%%%%%%%%%%%%%%%%%%%%%%%%%
\section{Evaluation using matrix functions}
\label{sec.matexp}
%%%%%%%%%%%%%%%%%%%%%%%%%%%%%%%%%%%%%%%%%%%%%%%%%%%%%%%%%%%%

We next describe different approaches for evaluating the approximation $X_t^m$
given in \eqref{eq:finite_sum} using matrix functions.

%%%%%%%%%%%%%%%%%%%%%%%%%%%%%%
\subsection{Diagonalization of $L$}
\label{sec.diagonalization}
%%%%%%%%%%%%%%%%%%%%%%%%%%%%%%
We first investigate the case where $L$ is diagonalizable,
i.e., $L = V D V ^{-1}$ where $D \in \mathbb{C}^{n \times n}$ is diagonal.
% As the functions $\varphi_{k,t}^{\cos}(z)$ defined
% by~\eqref{eq:phi_cos}
% are entire functions,
% they have converging series expansions around zero for all
% $z \in \mathbb{C}$,
% and consequently $\varphi_{k,t}^{\cos}(L) = V \varphi_{k,t}^{\cos}(D) V^{-1}$.
% Since $\varphi_{k,t}^{\cos}(D) = \mathrm{diag}(\varphi_{k,t}^{\cos}(d_{ii}))$
% this considerably lowers the computational complexity.
Then we can rewrite the
numerical approximation~\eqref{eq:phi_series} as
\begin{equation} \label{eq:diagonalized_method}
u_N(t) = \ee^{t L} u_0 + \sqrt{2} \, V \sum\limits_{k=1}^m \varphi_{k,t}^{\cos}(t D) V^{-1} \, Z_k.
\end{equation}
% Thus, if diagonalization of $L$ is feasible,
% and moreover $V^{-1}$ can be computed explicitly or is cheap to apply
% to a vector,
% \eqref{eq:diagonalized_method} gives a fast alternative.
% The diagonal matrices $\varphi_{k,t}^{\cos}(t D)$
% can be precomputed using the formula \eqref{eq:phi_cos} along with
% the vector $\ee^{t L} u_0$;
% whereupon repeatedly sampling with random $Z_k$ has an almost negligible cost.
The drawback of this approach is that not all matrices $L$
are diagonalizable and,
even if such a decomposition exists,
it will destroy any structure such as sparsity in $L$ and requires
large amounts of memory for larger matrices.
The computation is also rather expensive (around $25n^3$ flops).
%and does not enjoy the backward error properties of the
%Taylor series approach (section~\ref{sec.tayl-expans}).
However,
if $L$ is sufficiently small and diagonalizable it may be worth doing
the initial diagonalization to speed up the subsequent sampling.

An important special case is normal $L$. This means that $L$ is
unitarily diagonalizable, i.e., there exists a unitary
$V \in \mathbb{C}^{n \times n}$ such that $L = V D V^*$ for some
diagonal $D$.
If the elements of $Z_k \in \mathbb{R}^n$ are
i.i.d. normally distributed with variance 1,
then so are the elements
of the vector $V^* Z_k$ \cite[Thm.\;2.1.2.]{Tong}.
This implies that instead of \eqref{eq:diagonalized_method}
we may use the simplified sampling formula
\begin{equation} \label{eq:normal_method}
u_N(t) = \ee^{t L} u_0 + \sqrt{2} \,  V \left( \sum\limits_{k=1}^m \varphi_{k,t}^{\cos}(t D) \, \widetilde{Z}_k \right),
\end{equation}
where the elements of the vectors $\widetilde{Z}_k$ are
i.i.d. normally distributed with variance 1.

%%%%%%%%%%%%%%%%%%%%%%%%%%%%%%
\subsection{Linearization of the truncated SDE}
\label{sec.linear-trig}
%%%%%%%%%%%%%%%%%%%%%%%%%%%%%%
We next consider  methods for a general $L$. These methods can also exploit the possible sparsity of $L$.
We start by considering the truncated series \eqref{eq:truncated}
and denote the vector of the first $N$ basis functions by $y_N(t)$,
\begin{equation} \label{eq:Fourier_series}
  y_N(t) =
  \begin{bsmallmatrix}
  \cos( c_1 t) \\
  \vdots \\
  \cos(c_N t) \\
  \sin( c_1 t) \\
  \vdots \\
  \sin( c_N t)
\end{bsmallmatrix}.
\end{equation}

Using this notation it is clear that $y_N(t)$ satisfies the system of differential equations
\begin{equation*} % \label{eq:y_N_equation}
\dot{y}_N(t) = \begin{bmatrix}
           0 & C_N \\ -C_N & 0
          \end{bmatrix} y_N(t), \quad y_N(0) = \begin{bmatrix}  \mathbf{1} \\ 0 \end{bmatrix},
\end{equation*}
where $C_N =  \mathrm{diag}( c_1,\ldots, c_N)$ and $\mathbf{1}
= \begin{bmatrix} 1 & \ldots & 1 \end{bmatrix}^\mathrm{T}$.
This allows us to rewrite $y_N(t)$ as
%the solution to this ODE, i.e.
\begin{equation} \label{eq:y_N_exp}
y_N(t) = \exp \left( t \begin{bmatrix}
           0 & C_N \\ -C_N & 0
          \end{bmatrix} \right) \begin{bmatrix} \mathbf{1} \\ 0 \end{bmatrix}.
\end{equation}

Each realization of the SDE~\eqref{eq:semilinearSDE}
is simply a special case of the deterministic equation
analyzed in section~\ref{sec.deterministic}.
Indeed, we see that the solution presented in
Theorem~\ref{thm.sde_soln} is a special case of~\eqref{eq:phi_series}.
In order to compute our approximation to the solution efficiently
we must be able to evaluate matrix--vector products of
the form $\varphi_{k,t}^{\cos}(A)v$ and $\varphi_{k,t}^{\sin}(A)v$.
We reformulate this solution in terms of the matrix exponential.

%Then one can also take advantage of the very efficient routines
%available for its computation.

To begin our reformulation,
let us define $A_N,B_N \in \mathbb{R}^{n \times N}$ by
\begin{equation*}
A_N = \begin{bmatrix} a_1 & \ldots & a_N \end{bmatrix}  \quad \textrm{and} \quad
B_N = \begin{bmatrix} b_1 & \ldots & b_N \end{bmatrix}.
\end{equation*}
When $y_N(t)$ is defined as in \eqref{eq:Fourier_series}, we see that
\begin{align*}
  g_N(t)
  &= \begin{bmatrix} A_N & B_N \end{bmatrix} y_N(t) \\
  &= \begin{bmatrix} A_N & B_N \end{bmatrix}
\exp \left( t \begin{bmatrix}
           0 & -C_N \\ C_N & 0
         \end{bmatrix} \right) \begin{bmatrix} \mathbf{1} \\ 0  \end{bmatrix}
\end{align*}
Using this, we obtain the following result
with proof left to the appendix
(see also Thm. 2.1 in \cite{alhi11}).\\

%\begin{thm}[Representation of the solution using the matrix exponential]
\begin{thm} \label{thm:augmented}
  Let $L \in \mathbb{R}^{n \times n}$,
  $g_N(t)$ be the partial Fourier series defined by the
  coefficients $a_k,b_k,c_k$,
  and let $u_N(t)$ be the solution of \eqref{eq:semilinearODE_N}.
  Then,
  \begin{equation} \label{eq:exp_augmented}
    u_N(t) =
    \begin{bmatrix} I_n & 0 \end{bmatrix}
    \exp\left( t
      \begin{bsmallmatrix}
        L & A_N & B_N \\
        0 & 0 & -C_N \\
        0 & C_N  & 0
      \end{bsmallmatrix} \right)
    \begin{bsmallmatrix} u_0 \\  \mathbf{1} \\ 0 \end{bsmallmatrix}
  \end{equation}
  where $\mathbf{1} = \begin{bmatrix} 1 & \ldots & 1 \end{bmatrix}^\mathrm{T}$.
\end{thm}

For a more specific example,
let us consider the stochastic differential equation
\begin{equation*}
  \dot{u}_N(t) = Lu_N(t) + dW_t^N,
\end{equation*}
where $dW_t^N$ denotes the truncated Karhunen--Lo\'{e}ve expansion
of the Weiner process~\eqref{eq:weiner_khl}.
This is a special case of Theorem~\ref{thm.sde_soln} where $B = I_n$.
For this particular problem,
if $Z_N \in \mathbb{R}^{n \times N}$ is a matrix with
independent elements drawn from a normal $\mathcal{N}(0,1)$ distribution,
we would have
$A_N = [0,\dots,0]$,
$B_N = \sqrt{2} Z_N$, and
$C_N = \mathrm{diag} \left(  \frac{1}{2} \pi,\ldots,   \left(N -
    \frac{1}{2}\right) \pi  \right)$.

Since we have now simplified the computation of the
solution into the product of the matrix exponential
multiplied by a vector one can take advantage of
many efficient methods for its computation.
As an example we mention the scaling and squaring method~\cite{alhi11}
and the Krylov subspace methods~\cite{Saad}.
The best method for computing $u_N(t)$ for any
particular problem will depend largely upon the matrix $L$
(i.e. whether $L$ is small and dense or large and sparse)
and the accuracy required in the final solution.

We next consider a specific approach which exploits
the fact that only the matrix $B_N$ changes when using the expression
\eqref{eq:exp_augmented} to evaluate samples of $X_t^m$.

% We consider computing $u_N(t)$ via
% \begin{itemize}
% \item the explicit diagonalization of $L$,
% \item the matrix exponential with scaling and squaring.
% %\item the matrix exponential with Krylov subspace methods, and
% \end{itemize}

%
%
%%%%%%%%%%%%%%%%%%%%%%%%%%%%%%%%%%%%%%%%%%%%%%%%%%%%%%%%%%%%
\subsubsection{Sylvester equation approach}
\label{sec.sylvester}
%%%%%%%%%%%%%%%%%%%%%%%%%%%%%%%%%%%%%%%%%%%%%%%%%%%%%%%%%%%%

Since only the matrix $B_N$ changes for each realization of $X_t^m$,
it is only the $(1,2)$-block of size $n \times 2N$ in the exponential \eqref{eq:exp_augmented} that changes for each $X_t^m$.
We have the following result which can be used to efficiently compute
the $(1,2)$-block.
The proof is given in the appendix.
\begin{lem} % \label{lem:12block}

Let $L \in \mathbb{R}^{n \times n}$, $B \in \mathbb{R}^{n \times m}$ and $C \in \mathbb{R}^{m \times m}$. Then,
$$
\exp\left( t \begin{bmatrix} L & B \\ 0 & C \end{bmatrix} \right) = \begin{bmatrix} \ee^{tL} & X(t) \\ 0 & \ee^{tC} \end{bmatrix},
$$
where $X(t)$ satisfies the Sylvester equation
\begin{equation} \label{eq:sylvester}
	L X(t) - X(t) C = \ee^{tL} B - B \ee^{tC}.
\end{equation}
\end{lem}

When sampling, we set $C = \begin{bsmallmatrix} 0 & -C_N \\ C_N & 0 \end{bsmallmatrix}$ and $B = \begin{bsmallmatrix} 0 & B_N \end{bsmallmatrix}$.
Notice that $C$ has purely imaginary eigenvalues, so that if $A$ has its spectrum on the left half-plane, for example,
the spectra of $A$ and $C$ are well separated and the Sylvester equation \eqref{eq:sylvester} always has a solution.
This strategy has also been mentioned in~\cite[p.248]{high:FM}.

%
% (This strategy also recommended in N.~Higham's book, page 248:
% " If the spectra of A11 and A22 are well separated then it is best to compute eA11 and eA22 individually and obtain
% F12 from the block Parlett recurrence by solving a Sylvester equation of the form (9.7).", perhaps worth mentioning)

%%%%%%%%%%%%%%%%%%%%%%%%%%%%%%
\section{Numerical experiments}
\label{sec.numerical-experiments}
%%%%%%%%%%%%%%%%%%%%%%%%%%%%%%
We are now ready to test our novel methods against
the Euler--Maruyama and backward Euler--Maruyama scheme.
We use two illustrative stochastic differential equations
to compare the methods.
%with varying levels of complexity.
The first one is a small dimensional equation so that the diagonalization approach described in Sec.~\ref{sec.diagonalization}
can be used. Moreover, the equation is non-stiff so we compare it to the explicit Euler--Maruyama scheme.
The second SDE is a large-dimensional stiff equation, and we compare the Sylvester equation based approach of
Sec.~\ref{sec.sylvester} to the backward Euler--Maruyama scheme.
As a metric for comparison we use the weak convergence of the second order moment $\mathbb{E} \norm{X_t}^2 $.

All experiments in this section were performed
on a laptop machine with a Intel Core i5 (3.1~GHz)
with 16GB of RAM.
Computations were performed with MATLAB~2016b.

% \textbf{Some brief description of the experiments and metrics
%   used for comparison here}.

The sampling errors of the computed quantities decay like $1/\sqrt{N}$, where $N$ is the size of the sample.
For details, see~\cite[Sec.\;1.9]{kloeden}. Thus, when comparing the convergences of different methods,
attention has to be paid to the selection of large enough $N$.

The MATLAB code for the experiments is provided in the supplementary material.

%%%%%%%%%%%%%%%%%%%%%%%%%%%%%%
\subsection{Turbulent diffusion}
\label{sec.turbulent-diffusion}
%%%%%%%%%%%%%%%%%%%%%%%%%%%%%%
For our first example we consider the following small example of
turbulent diffusion taken from Kloeden and Platen~\cite[Sec.~7]{kloeden}.
% Let $X_t \in \mathbb{R}^3$ denote the position of a fluid particle at time $t$ and $V_t$ its velocity.
% One of the early models for such turbulent diffusion due to Obukhov was
% \begin{equation*}
%   dX_t = V_tdt,\quad dV_t = -\frac{1}{T}V_tdt + \sigma dW_t.
% \end{equation*}
% Clearly here $V_t$ is an Ornstein--Uhlenbeck process. %The main criticism of this model is that the large and small scale
% %effects happen simultaneously which are not well accounted for in this simple model.
% Atmospheric scientists have since proposed more complicated systems to account for %this
% multi-scale behaviour.
For variables $V_t^{(1)},V_t^{(1)} \in \mathbb{R}^3$,
the equations describing the system are given by
%proposed the following model which will be used throughout this experiment.
\begin{equation*} %\label{eq:diffusion_motion}
\begin{aligned}
dV_t^{(1)} &= \left( - \frac{1}{T_1} V_t^{(1)} - \beta( V_t^{(1)} - V_t^{(2)} ) \right) dt + \sigma_1 dW_t^{(1)} \\
dV_t^{(2)} &= \left( - \frac{1}{T_2} V_t^{(2)} + \beta( V_t^{(1)} - V_t^{(2)} ) \right) dt + \sigma_2 dW_t^{(2)},
\end{aligned}
\end{equation*}
where $\sigma_1,\sigma_2,T_1,T_2$ are constants that determine the behaviour of the system.
This 6 dimensional system can be reformulated into a single equation: if we denote $V = [V^{(1)},~V^{(2)}]^T$ then
\begin{equation} \label{eq:diffusion_motion}
  dV_t =
  \begin{bsmallmatrix}
    -(\frac{1}{T_1} + \beta) I & \beta I\\
    \beta I & -(\frac{1}{T_2} + \beta)I
  \end{bsmallmatrix} V_t +
  \begin{bsmallmatrix}
    \sigma_1 I & 0 \\
    0 & \sigma_2 I
  \end{bsmallmatrix} dW_t.
\end{equation}

We observe that the system is of the
%that~\eqref{eq:diffusion_motion} is of the
form~\eqref{eq:semilinearSDE} with a symmetric coefficient matrix $L$.
This allows us to apply the efficient sampling formula~\eqref{eq:normal_method}
based on the diagonalization of normal matrices
(see section~\ref{sec.diagonalization}).

Our aim is to compare this diagonalization procedure
based upon the Karhunen--Lo\'{e}ve expansion
to the Euler--Maruyama time-stepping scheme
$$
V_{i+1} = V_i + \Delta t L  V_i + B \Delta W_i,
$$
where $\Delta W_i$'s are independent normally distributed random variables with covariance $\Delta t I_n$.
The coefficient matrices $L$ and $B$ are given in \eqref{eq:diffusion_motion}.

\begin{figure}[t!]
 \begin{center}
 \includegraphics[scale=0.45]{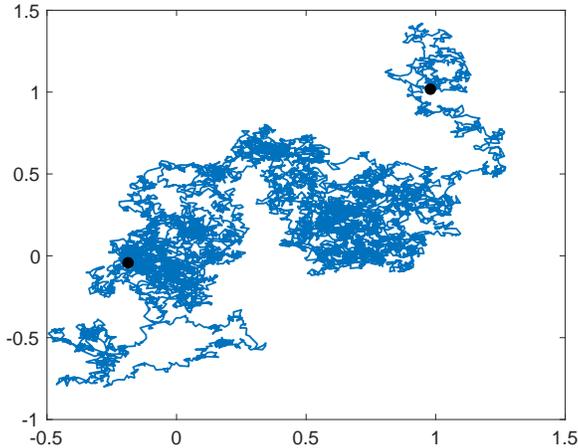}
 \end{center}
 \caption{First two components of a single trajectory
   $V_t^m$ with $m=1000$ time discretization points.
   The black dots depict the initial and end points,
   starting at $[1~1]^T$.}
 \label{fig:trajectories}
\end{figure}

For this experiment we set the parameters $T_1=T_2 = 0.5$,
$\sigma_1 = \sigma_2 = 1$, and $\beta = 2$.
The initial value is set to
$V_0 = \begin{bmatrix} 1 & \ldots & 1 \end{bmatrix}^T$.
To give some insight into the typical behaviour a fluid particle might
have under these conditions,
we have plotted a single particle trajectory
(projected onto a 2D plane) in Figure~\ref{fig:trajectories}.
This trajectory was computed using the Euler--Maruyama scheme with
$m = 1000$ time discretization points.

\begin{figure}[t!]
 \begin{center}
 \includegraphics[scale=0.45]{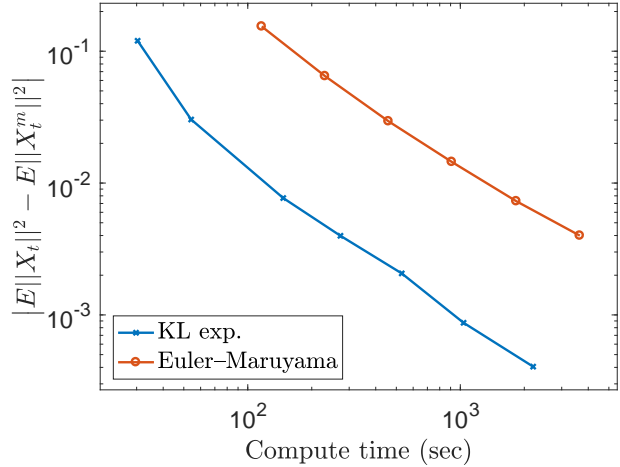}
 \end{center}
 \caption{Compute times vs. the weak error $\abs{\mathbb{E}(\norm{X_t}^2) - \mathbb{E}(\norm{X_t^m}^2)}$ for the KL-expansion based method
 for different lengths of the expansion and Euler--Maruyama method
 with different number of time steps. Each estimate of $\mathbb{E}(\norm{X_t^m}^2)$ is an average of
 $10^7$ samples.}
 \label{fig:cpu1}
\end{figure}

We denote by $m$ both the number of
time discretization points in the Euler--Maruyama method and the
length of the Karhunen--Lo\'{e}ve expansion.
We evaluate the KL expansion based method for values $m=10,40,160,320,640,1280,2560$.
The Euler--Maruyama method is evaluated for $m=10,20,40,80,160,320$.
In Figure~\ref{fig:cpu1} we denote by $X_t^m$ the approximation of both methods.

%%%%%%%%%%%%%%%%%%%%%%%%%%%%%%
%%%%%%%%%%%%%%%%%%%%%%%%%%%%%%

%%%%%%%%%%%%%%%%%%%%%%%%%%%%%%
%%%%%%%%%%%%%%%%%%%%%%%%%%%%%%

For each estimate of $\mathbb{E}(\norm{X_t^m}^2)$ we draw $10^7$ samples and
the reference value $\mathbb{E}(\norm{X_t}^2)$ is computed in high precision using the expression given in Corollary~\ref{cor:normal}.
There is little difference between the two approaches in terms of their convergence,
but the Karhunen--Lo\'{e}ve approach (using the sampling method described in section~\ref{sec.diagonalization})
is much more efficient: both runtimes scale linearly with the sample size but
the Karhunen--Lo\'{e}ve approach is around 10 times faster.

%%%%%%%%%%%%%%%%%%%%%%%%%%%%%%
 \subsection{A finite difference discretization of a
   heat equation with additive noise}
 \label{sec.advection-diffu}
\begin{figure}[t!]
 \begin{center}
 \includegraphics[scale=0.45]{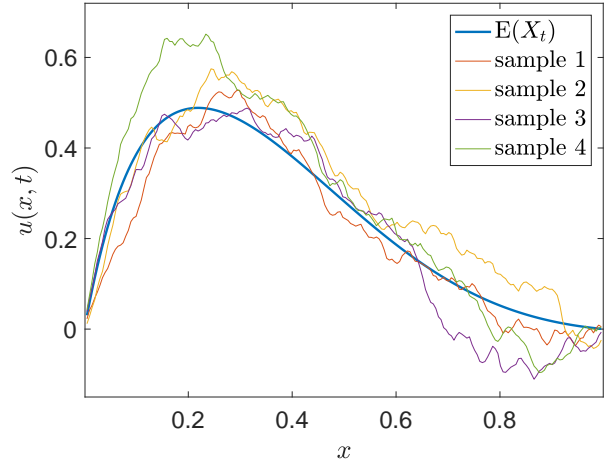}
 \end{center}
 \caption{$\mathbb{E}[X_t] $ and 4 realizations of $X_t$ at $t=0.4$. }
 \label{fig:large_ends}
\end{figure}

%%%%%%%%%%%%%%%%%%%%%%%%%%%%%%
 In our second experiment we consider a
% finite difference discretization of the
 one dimensional stochastic partial differential equation
\begin{equation*} % \label{eq:adv_diff}
\begin{aligned}
\frac{\partial}{\partial t} y(x,t) &= \varepsilon \frac{\partial^2}{\partial x^2} y(x,t)+ \alpha \frac{\partial}{\partial x} y(x,t)+
\beta \frac{\partial^2 W}{\partial t \partial x}, \\
 y(0,x) &= y_0(x), \quad 0 \leq x \leq 1, \\
 y(t,0) &= y(t,1) = 0, \quad t \geq 0.
\end{aligned}
\end{equation*}
Here $\frac{\partial^2 W}{\partial t \partial x}$ denotes mixed second order derivative of the so called Brownian sheet~\cite{anz98}.
The initial value $y_0(x)$ is given by
\begin{equation*}
	y_0(x)=\begin{cases}
		2x, &\text{ if } 0 \leq x \leq \frac{1}{2} \\
		2-2x, &\text{ if } \frac{1}{2} \leq x \leq 1.
	\end{cases}
\end{equation*}
Spatial discretization using finite differences leads to the SDE (see e.g.~\cite{anz98})
\begin{equation} \label{eq:pdesde}
	\begin{aligned}
		d X =  L  X \; dt+ \frac{\beta}{\sqrt{\Delta} x} dW_t, \quad X(0) = X_0,
	\end{aligned}
\end{equation}
where $\Delta x = L/(n+1)$, $L = \varepsilon \Delta_n + \alpha \nabla_n$,
\begin{equation*}
\Delta_n = \tfrac{1}{(\Delta x)^2} \begin{bsmallmatrix}
    -2 & 1 & & \\
     1 & \ddots & \ddots & \\
      & \ddots & \ddots & 1 \\
      & & 1 & -2 \end{bsmallmatrix}, \,
\nabla_n = \tfrac{1}{2 \Delta x} \begin{bsmallmatrix}
      & 1 & & \\
     -1 &  & \ddots & \\
      & \ddots &   & 1 \\
      & & -1 &  \end{bsmallmatrix},
\end{equation*}
and $X_0 \in \mathbb{R}^n$ is the discretization of the initial value $y_0(x)$.
We set $n=200$, $\veps=0.1$, $\alpha=-1.0$, $\beta=0.1$, and integrate up to $t=0.4$.
Figure~\ref{fig:large_ends} depicts 4 random samples and the expectation at $t=0.4$.

%%%%%%%%%%%%%%%%%%%%%%%%%%%%%%
%%%%%%%%%%%%%%%%%%%%%%%%%%%%%%

We evaluate each realization of the truncated KL-expansion using the augmented exponential and
Sylvester equation technique described in Section~\ref{sec.linear-trig}.
This means that at each step we solve a Sylvester equation of the form \eqref{eq:sylvester}.
As the coefficient matrix is sparse, the Sylvester equation is
efficiently solved by vectorizing the equation
and reusing the sparse LU factors throughout the sampling process.

As the SDE \eqref{eq:pdesde} is now stiff (see~\cite{kloeden}), we compare it
%based upon the Karhunen--Lo\'{e}ve expansion
to the backward Euler--Maruyama method
\begin{equation} \label{eq:bem}
	X_{i+1} = (I - \Delta t L)^{-1} X_i + B \Delta W_i,
\end{equation}
where $\Delta W_i$'s are independent normally distributed random vectors with covariance $\Delta t I_n$
and $B = \frac{\beta}{\sqrt{\Delta} x} I_n$.

In the implementation of the backward Euler--Maruyama method, we use a precomputed sparse LU-factorization for $I - \Delta t L$ to evaluate the time steps.

\begin{figure}[t!]
 \begin{center}
 \includegraphics[scale=0.45]{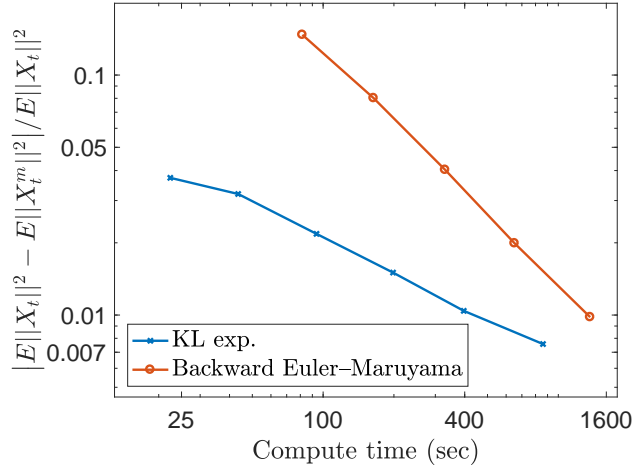}
 \end{center}
 \caption{Compute times vs. the relative weak error $\tfrac{\abs{\mathbb{E}(\norm{X_t}^2) - \mathbb{E}(\norm{X_t^m}^2)}}{ \mathbb{E}(\norm{X_t}^2) }$
 for the KL-expansion based method and backward Euler--Maruyama method
 with different number of time steps. Each estimate of $\mathbb{E}(\norm{X_t^m}^2)$ is an average of
 $10^5$ samples.}
 \label{fig:cpu2}
\end{figure}

We evaluate the KL expansion based method for the expansion length $m=1,2,\ldots,2^5$.
The backward Euler--Maruyama method is evaluated for number of time steps $m=50,100,\ldots,800$.
In Figure~\ref{fig:cpu2} we denote by $X_t^m$ the approximation of both methods.

The Karhunen--Lo\'{e}ve approach (using the Sylvester equation approach described in Sec.~\ref{sec.sylvester})
is again more efficient (see Figure~\ref{fig:cpu2}). However, although both methods have the weak order of convergence $O(\frac{1}{m})$,
the compute time of KL approach grows quadratically with $m$ which explains the difference in the slopes of the lines in Figure~\ref{fig:cpu2}.
Here improvements could be made, e.g., for evaluating the right hand side of the Sylvester equation \eqref{eq:sylvester}.

%%%%%%%%%%%%%%%%%%%%%%%%%%%%%%
%%%%%%%%%%%%%%%%%%%%%%%%%%%%%%

%%%%%%%%%%%%%%%%%%%%%%%%%%%%%%
\section{Conclusions}
\label{sec:conclusion}
%%%%%%%%%%%%%%%%%%%%%%%%%%%%%%

We have proposed a novel approach for sampling linear SDEs which exploits
matrix functions and efficient numerical linear algebraic subroutines.
Moreover, we have provided a convergence theory for the method which
shows both the strong and weak convergence speeds. In numerical examples we showed that the method is very competitive
both for non-stiff and stiff systems of equations, when comparing
against the Euler--Maruyama and backward Euler--Maruyama method, respectively.
As future work, we are interested in optimizing the various linear algebraic
subproblems using Krylov subspaces, for example.
Furthermore, our approach should be generalizable to
inhomogeneous linear SDEs (see e.g.~\cite[Sec.\;6]{saso19}).

\clearpage
\newpage

%%%%%%%%%%%%%%%%%%%%%%%%%%%%%%
\bibliographystyle{plain}
\bibliography{strings,manuscript}

\appendix

\section{Proofs}

\begin{customthm}{4}
We have
\begin{enumerate}%[label=(\roman*)]

\item
$
\mathbb{E}(X_t) = \ee^{t L} X_0.
$
\item
$
\mathbb{E} \left( \norm{X_t}^2 \right) = \norm{\ee^{t L} X_0}^2 + 2 \sum\limits_{k=1}^\infty \norm{ \varphi_{k,t}^{\cos}(L) B }_F^2,
$
where $\norm{ \cdot }_F$ is the Frobenius norm.
\end{enumerate}
\end{customthm}
\begin{proof}

We see from \eqref{eq:truncated_integral} that for all $m>0$
$$
\mathbb{E}\left( \int_0^t \ee^{(t-s) L} B \, d W_s^m \right) = \sqrt{2} \sum_{k=1}^m  \varphi_{k,t}^{\cos}(L) B \, \mathbb{E}\left( Z_k \right) = 0 \quad
$$
which implies the first claim.

Since the elements of $Z_i$'s are i.i.d. $\mathcal{N}(0,1)$ - distributed, it is easily verified that
for all $A \in \mathbb{R}^{n \times n}$ we have $\mathbb{E}(Z_i^T A Z_j) = \delta_{ij} \, \mathrm{tr}(A)$,
where $\mathrm{tr(A)}$ denotes the trace
of $A$. Thus, it follows from \eqref{eq:series_representation} that
\begin{equation*}
\begin{aligned}
 \mathbb{E} \norm{X_t}^2 & = \norm{ \ee^{t L} X_0 }^2 + 2 \, \sum_{k=1}^\infty
 \mathbb{E} \left( \left(\varphi_{k,t}^{\cos}(L) B Z_k\right)^T \varphi_{k,t}^{\cos}(L) B Z_k \right) \\
 & = \norm{ \ee^{t L} X_0 }^2  + 2  \, \sum_{k=1}^\infty
 \mathrm{tr} \left( \left(\varphi_{k,t}^{\cos}(L) B \right)^T \varphi_{k,t}^{\cos}(L) B \right).
 \end{aligned}
\end{equation*}
 The second claim follows then the fact that for all $A \in \mathbb{R}^{n \times n}$, $\mathrm{tr}(A^T A) = \norm{A}_F^2$.
\end{proof}

\begin{customthm}{6}
 Let the linear operator $L \in \mathbb{R}^{n \times n}$
 in \eqref{eq:semilinearSDE} be negative semidefinite.
 The error introduced by the approximation~\eqref{eq:finite_sum}
 with $m$ terms satisfies
 \begin{equation*}
	 \begin{aligned}
	    \mathbb{E}\left(\norm{X_t - X^m_t}^2 \right) &= 2 \sum_{k={m+1}}^\infty \normm{\varphi_{k,t}^{\cos}(L) B }_F^2 \\
		& \leq \frac{2 \norm{B}^2}{\pi^2} \frac{n}{m-1}.
	 \end{aligned}
 \end{equation*}
\end{customthm}
\begin{proof}
From the representation \eqref{eq:series_representation} and \eqref{eq:finite_sum}, it follows that
 %
%Since $\mathbb{E}(Z_i^T A Z_j) = 0$ for all $i\neq j$ and for all $A \in \mathbb{R}^{n \times n}$,
  \begin{equation*} % \label{eq:exp_ineq}
  \mathbb{E}\big(\normm{X_t - X^m_t}^2 \big) =
  \mathbb{E}\left(\normm{ \sqrt{2} \sum_{k={m+1}}^\infty  \varphi_{k,t}^{\cos}(L) B \, Z_k }^2  \right).
  \end{equation*}
Following the lines of the proof of Theorem~\ref{thm:variance_thm}, we see
 \begin{equation*}% \label{eq:exp_ineq}
  \mathbb{E}\big(\normm{X_t - X^m_t}^2 \big) = 2 \sum_{k={m+1}}^\infty \normm{\varphi_{k,t}^{\cos}(L) B }_F^2.
  \end{equation*}
% Furthermore, for all matrices $A \in \mathbb{R}^{n \times n}$ and all $Z_k$, $k \geq 1$, we have
%$$
%\mathbb{E}\norm{A Z_k}^2 \leq \mathbb{E} \norm{A}^2 \norm{Z_k}^2 = n \norm{A}^2.
%$$
%Thus, from \eqref{eq:exp_ineq}
For all $A\in \mathbb{R}^{n \times n}$ it holds that $\norm{A}_F \leq \sqrt{n} \norm{A}_2$, and therefore
$$
\mathbb{E}\big(\normm{X_t - X^m_t}^2 \big) \leq 2 n \norm{B}^2 \sum_{k={m+1}}^\infty \normm{\varphi_{k,t}^{\cos}(L)}^2.
$$
Using Lemma~\ref{lem:norm_bound} to bound $\normm{\varphi_{k,t}^{\cos}(L)}$ gives
\begin{equation*}
\begin{aligned}
%\mathbb{E}\big(\normm{X_t - X^m_t} \big)^2
\sum_{k={m+1}}^\infty \normm{\varphi_{k,t}^{\cos}(L)}^2 & \leq  \sum_{k={m+1}}^\infty \frac{1}{\big((k-\frac{1}{2})\pi\big)^2} \\
& \leq   \sum_{k={m}}^\infty \frac{1}{\big(k \pi\big)^2} \\
&\leq  \int\limits_{m-1}^\infty \frac{1}{(\pi x)^2} \, \dd x \\
&  = \frac{1}{\pi^2 (m-1)}.
\end{aligned}
\end{equation*}
 \end{proof}

%\begin{thm}[Representation of the solution using the matrix exponential]
\begin{customthm}{10}
  Let $L \in \mathbb{R}^{n \times n}$,
  $g_N(t)$ be the partial Fourier series defined by the
  coefficients $a_k,b_k,c_k$,
  and let $u_N(t)$ be the solution of \eqref{eq:semilinearODE_N}.
  Then,
  \begin{equation*}
    u_N(t) =
    \begin{bmatrix} I_n & 0 \end{bmatrix}
    \exp\left( t
      \begin{bmatrix}
        L & A_N & B_N \\
        0 & 0 & -C_N \\
        0 & C_N  & 0
      \end{bmatrix} \right)
    \begin{bmatrix} u_0 \\  \mathbf{1} \\ 0 \end{bmatrix}
  \end{equation*}
  where $\mathbf{1} = \begin{bmatrix} 1 & \ldots & 1 \end{bmatrix}^\mathrm{T}$.
\end{customthm}

\begin{proof}
  The claim follows from the fact that for any square matrices
  $X_1$ and $X_2$~\cite[pp. 248]{high:FM}
\begin{equation*}
  \exp \left( t \begin{bmatrix} X_1 & X_3 \\ 0 & X_2 \end{bmatrix} \right) = \begin{bmatrix} \ee^{tX_1} & \int\limits_0^t \ee^{(t-s)X_1} X_3 \, \ee^{s X_2} \, \dd s \\ 0 & \ee^{t X_2} \end{bmatrix}.
\end{equation*}
We can select $X_1 = L$, $X_2 =
\smatrix{0 & -C_N \\ C_N & 0}
$,
and $X_3 = [A_N, B_N]$
before combining this result with the substitution
\begin{equation*}
  \exp \left( t
    \begin{bmatrix} 0 & -C_N \\ C_N & 0 \end{bmatrix} \right) =
  \begin{bmatrix} \cos( t \, C_N) & -\sin(t \, C_N) \\
\sin(t \, C_N) & \cos(t \, C_N) \end{bmatrix}.
\end{equation*}
\end{proof}

\begin{customlemma}{11}
Let $L \in \mathbb{R}^{n \times n}$, $B \in \mathbb{R}^{n \times m}$ and $C \in \mathbb{R}^{m \times m}$. Then,
$$
\exp\left( t \begin{bmatrix} L & B \\ 0 & C \end{bmatrix} \right) = \begin{bmatrix} \ee^{tL} & X(t) \\ 0 & \ee^{tC} \end{bmatrix},
$$
where $X(t)$ satisfies the Sylvester equation
\begin{equation} \label{eq:sylvester}
	L X(t) - X(t) C = \ee^{tL} B - B \ee^{tC}.
\end{equation}
\end{customlemma}
\begin{proof}
Since every matrix commutes with its exponential, it holds
\begin{equation*}
\begin{bmatrix} L & B \\ 0 & C \end{bmatrix} 	\begin{bmatrix} \ee^{tL} & X(t) \\ 0 & \ee^{tC} \end{bmatrix} =
	\begin{bmatrix} \ee^{tL} & X(t) \\ 0 & \ee^{tC} \end{bmatrix} \begin{bmatrix} L & B \\ 0 & C \end{bmatrix}.
\end{equation*}
The $(1,2)$-block of this matrix equation gives the Sylvester equation \eqref{eq:sylvester}.
\end{proof}

\end{document}